\documentclass[12pt]{article}

\usepackage{pstricks,pst-node}
\usepackage{amsmath,amsthm}
\usepackage{amssymb}
\usepackage{color}
\usepackage{xspace}
\usepackage[colorlinks=true,
linkcolor=green,
filecolor=brown,
citecolor=green]{hyperref}

\def\eps{\epsilon}
\def\a{\ensuremath{\mathcal A}\xspace}
\def\b{\ensuremath{\mathcal B}\xspace}
\def\c{\ensuremath{\mathcal C}\xspace}

\def\n{N\xspace}
\def\gf{generating function\xspace}
\def\gl{ground level\xspace}

\begin{document}
\newtheorem*{theorem}{Theorem} 

\begin{center}
{\Large
A Note on Generalized Narayana Numbers                           \\ 
}
\vspace{2mm}
David Callan  
\vspace{2mm}

\end{center}

\begin{abstract}
We give a simple proof that the the number of Dyck paths of semilength $n$ with $i$ returns to \gl and $j$ peaks is the generalized Narayana number $\frac{i}{n} \binom{n}{j} \binom{n - i - 1}{j - i}$.
\end{abstract}

The Narayana numbers, $\frac{1}{n}\binom{n}{k}\binom{n}{k-1},\ 1\le k \le n$, form sequence \htmladdnormallink{A001263}{http://oeis.org/A001263} in OEIS \cite{oeis}, and count Dyck paths \cite[page 7]{cat2015} by number of peaks. 
Various generalizations have been considered \cite{sulanke2004,barry2011,kruchinin2021}. Here, we consider the 3-index array $\n_{i}(n,k)$ defined by 
\[
\n_{i}(n,k) =  \frac{i}{n} \binom{n}{j} \binom{n - i - 1}{j - i}
\]
for $ 1 \le i \le j \le n$.
Below is a table of values of $\n_i(n,j)$ for small $i,n ,j$.

 \[
\begin{array}{c|ccccc}
	n^{\textstyle{\backslash \:j}} &  1 & 2 & 3 & 4 & 5 \\
\hline 
	1&    1 &   &   &  &  \\
	2&    1 & 0 &   &   & \\
 	3&    1 & 1 & 0  &  &  \\
	4&    1 & 3 & 1 &  0 & \\
	5&    1 & 6 & 6 &  1 & 0\\ 
	\multicolumn{5}{c}{\raisebox{-1ex}{\quad $i=1$}}
 \end{array}
 \hspace*{20mm}
 \begin{array}{c|ccccc}
	n^{\textstyle{\backslash \: j}} &2 & 3 & 4 & 5 & 6 \\
\hline 
	2&    1 &     & & & \\
 	3&    2 & 0    & &  & \\
	4&    3 & 2 & 0&  &\\
	5&    4 & 8 & 2 & 0 & \\ 
	6&    5 & 20 & 15 & 2 & 0\\ 
	\multicolumn{5}{c}{\raisebox{-1.5ex}{\quad $i=2$}}
 \end{array}
\]

\vspace*{5mm}
 \[
\begin{array}{c|ccccc}
	n^{\textstyle{\backslash \:j}} & 3 & 4 & 5 & 6 & 7\\
\hline 
	3&    1 &   &   &  &  \\
	4&    3 & 0  &   &    &  \\
 	5&    6 & 3 & 0  &   & \\
	6&    10 & 15 & 3 & 0 & \\
	7&    15 & 45 & 27 &  3  & 0\\ 
	\multicolumn{5}{c}{\raisebox{-1ex}{$i=3$}}
 \end{array}
 \hspace*{20mm}
 \begin{array}{c|ccccc}
	n^{\textstyle{\backslash \:j}} & 4 & 5 & 6 & 7 & 8 \\
\hline 
	4&    1  &   & & & \\
 	5&    4  & 0 & & &  \\
	6&    10 & 4 & 0  &  &\\
	7&    20 & 24 & 4 & 0 & \\ 
	8&    35 & 84 & 42 & 4 & 0 \\ 
	\multicolumn{5}{c}{\raisebox{-1.5ex}{$i=4$}}
 \end{array}
\]

\vspace*{3mm}
\centerline{Values of $\n_i(n,j)$ for small $i,n,j$}
\vspace*{2mm}
The triangular array for $i=2$ is sequence \htmladdnormallink{A108838}{http://oeis.org/A108838} and the array for $i=3$ (with rows reversed)  is 
\htmladdnormallink{A281293}{http://oeis.org/A281293}.
\begin{theorem}
$\n_i(n,j)$ is the number of Dyck paths of semilength $n$ with $i$ returns to \gl and $j$ peaks.
\end{theorem}
\begin{proof} The proof is a nice application of an involution $\phi$ on Dyck paths due to Emeric Deutsch \cite{deutschInvol}, a bijection from Dyck paths to parallelogram polyominos \cite[Chapter 2, Item 57]{cat2015}, and the Lindstr\"{o}m-Gessel-Viennot theorem to count a set of nonintersecting lattice paths with given endpoints \cite[search Lindstrom-Gessel-Viennot lemma\,]{wiki}.

The involution $\phi$ is defined recursively using the first return decomposition $UP_1DP_2$ of a nonempty Dyck path  as follows (where $U=(1,1)$ is an upstep, $D=(1,-1)$ is a downstep, $P_1,\,P_2$ are Dyck paths, and $\eps$ denotes the empty path): 
\[ 
\phi(\eps)=\eps, \quad \phi(UP_1DP_2)=U\phi(P_2)D\phi(P_1)\,.
\]
The map $\phi$ sends number of returns to length of initial ascent, and if a Dyck path $P$ has semilength $n$ and $j$ peaks, then $\phi(P)$ has $n+1-j$ peaks (this explains the symmetry in the distribution of peaks in Dyck paths).

Let $\a_i(n,j)$ denote the set of paths being counted. Then $\phi:\,\a_i(n,j) \to 
\b_i(n,j)$ where $\b_i(n,j)$ is the set of Dyck paths $P$ of semilength $n$ with 
initial ascent of length $i$ and $n+1-j$ peaks. 
See below for an example with $n=10,\,i=4,\,j=6$.
\begin{center}
\begin{pspicture}(-1,-.5)(48,1)
\psset{unit=.3cm}
\psline[linewidth=.05cm](0,0)(1,1)(2,2)(3,3)(4,2)(5,1)(6,2)(7,1)(8,0)(9,1)(10,2)(11,1)(12,2)(13,3)(14,2)(15,1)(16,0)(17,1)(18,0)(19,1)(20,0)

\rput(10,-1.5){\small{Dyck path in $\a_i(n,j)$}}
\rput(10,-3){\small{with $n=10,\,i=4,\,j=6$}}
\rput(35,-1.5){\small{Dyck path in $\b_i(n,j)$}}
\rput(35,-3){\small{first ascent length is 4; 5 peaks}}
\rput(22.5,2){\textrm{\large $\xrightarrow{\phi}$}}

\rput(48,2){\textrm{\large $\rightarrow$}}  

\psline[linecolor=gray](0,0)(20,0)
\psline[linecolor=gray](25,0)(45,0)

\psset{dotsize=3pt 0}
\psdots(0,0)(1,1)(2,2)(3,3)(4,2)(5,1)(6,2)(7,1)(8,0)(9,1)(10,2)(11,1)(12,2)(13,3)(14,2)(15,1)(16,0)(17,1)(18,0)(19,1)(20,0)

\psline[linewidth=.05cm](25,0)(26,1)(27,2)(28,3)(29,4)(30,3)(31,2)(32,1)(33,2)(34,3)(35,2)(36,3)(37,2)(38,1)(39,0)(40,1)(41,2)(42,1)(43,0)(44,1)(45,0)

\psdots(25,0)(26,1)(27,2)(28,3)(29,4)(30,3)(31,2)(32,1)(33,2)(34,3)(35,2)(36,3)(37,2)(38,1)(39,0)(40,1)(41,2)(42,1)(43,0)(44,1)(45,0)

\end{pspicture}
\end{center}

\begin{center}
\begin{pspicture}(-1,-1.5)(5,7.5)
\psset{unit=1cm}
\psline[linewidth=.05cm](0,0)(1,0)(1,1)(1,2)(2,2)(3,2)(3,3)(3,4)(4,4)(4,5)(5,5)(5,6)
\psset{dotsize=3pt 0}
\psdots(0,0)(1,0)(1,1)(1,2)(2,2)(3,2)(3,3)(3,4)(4,4)(4,5)(5,5)(5,6)

\psline[linewidth=.05cm](0,0)(0,1)(0,2)(0,3)(0,4)(1,4)(1,5)(2,5)(3,5)(3,6)(4,6)(5,6)
\psdots(0,0)(0,1)(0,2)(0,3)(0,4)(1,4)(1,5)(2,5)(3,5)(3,6)(4,6)(5,6)

\rput(2.1,-0.7){\small{The parallelogram polyomino in $\c_i(n,j)$ determines }}
\rput(2.1,-1.2){\small{a pair of nonintersecting lattice paths, $A_1B_1$ and $A_2B_2$}}

\rput(1.5,3.7){\footnotesize{$A_1=(1,i)$}}
\rput(2,0){\footnotesize{$A_2=(1,0)$}}
\rput(4.8,6.3){\footnotesize{$B_1=(n-j,j)$}}
\rput(6.5,4.6){\footnotesize{$B_2=(n+1-j,j-1)$}}

\psset{dotsize=5pt 0}
\psdots(1,4)(1,0)(5,5)(4,6)
 
\end{pspicture}
\end{center}

The bijection from $\b_i(n,j)$ to parallelogram polyominos can be summarized as follows: Split the path into its ascents and descents. Consider each ascent in turn, left to right, as a vertical segment pointing north and rotate its topmost step clockwise $90^\circ$ so that it points east. Then concatenate and prepend a vertical step to get the upper boundary. Do likewise with the descents except rotate the bottom step, and append a vertical step  to get the lower boundary. The resulting set $\c_i(n,j)$ of parallelogram polyominos is characterized as those whose upper right vertex is $(n+1-j,j)$ and whose initial vertical segment on the upper boundary has length $i$. By deleting its mandatory steps, such a parallelogram polyomino is equivalent to a pair of nonintersecting lattice paths of unit vertical and horizontal steps, the upper path running from $A_1=(1,i)$ to $B_1=(n-j,j)$ and the lower path running from $A_2=(1,0)$ to $B_2=(n+1-j,j-1)$. Now Lindstr\"{o}m-Gessel-Viennot says the number of such pairs is 
\vspace*{-2mm}
\[
\#\, (\textrm{paths }A_1 \to B_1) \ \times\ \#\, (\textrm{paths }A_2 \to B_2)\ - \ 
\#\, (\textrm{paths }A_1 \to B_2) \ \times\ \#\, (\textrm{paths }A_2 \to B_1)
\]
\[
= \binom{n-j-1}{j-i} \times \binom{n-j}{j-1} - \binom{n-j}{j-1-i} \times \binom{n-j-1}{j}\  = \ \frac{i}{n} \binom{n}{j}   \binom{n-i-1}{j-i}
\]
\end{proof}
It is routine to show the \gf for these generalized Narayana numbers is given by 
\[
1+ \sum_{n\ge j\ge i \ge 1} \n_i(n,j) x^n y^i z^j =
\frac{2}{2-y \left(1-x (1-z)-\sqrt{1-2 x (1+z) + x^2 (1-z)^2} \right)}\,.
\]

\vspace*{8mm}

\centerline{\large{\textbf{Acknowledgment}}} 
The connection between the generalized Narayana numbers $\n_i(n,j)$ and the paths $\b_i(n,j)$ above was first observed by Werner Schulte \cite{schulte}.

\textsc{Department of Statistics, University of Wisconsin-Madison}

\end{document}